\documentclass{article}
\usepackage{arxiv}
\usepackage{graphicx} 
\usepackage{cmap}			 
\usepackage{mathtext} 		 
\usepackage[T2A]{fontenc}	 
\usepackage[utf8]{inputenc}	 
\usepackage[russian]{babel}	 
\usepackage{amsmath,amsfonts,amssymb,amsthm,mathtools}
\usepackage{icomma}
\usepackage{mdframed}
\usepackage{lipsum}
\usepackage{indentfirst}
\usepackage{url}
\usepackage{hyperref}

\newmdtheoremenv{theo}{Теорема}

\renewcommand{\leq}{\leqslant}
\renewcommand{\geq}{\geqslant}

\renewcommand{\ge}{\geqslant}

\date{}

\begin{document}

\title{Правила остановки градиентного метода для седловых задач с двусторонним  условием Поляка-Лоясиевича}

\author{
    Муратиди Александр Янисович \\
    МФТИ \\
    Москва, Россия \\
    \texttt{muratidi.aia@phystech.edu} \\
    \And
    Стонякин Федор Сергеевич \\
    МФТИ \\
    Москва, Россия \\
    КФУ им. В.\,И.\,Вернадского \\
    Симферополь, Россия \\
    \texttt{fedyor@mail.ru}
}

\maketitle


\begin{abstract}
Статья посвящена некоторым вопросам, связанным с численными методами для седловых задач с двусторонним условием градиентного доминирования Поляка--Лоясиевича \cite{yang}. Известно, что на классе достаточно гладких минимизационных задач с двусторонним условием градиентного доминирования градиентный метод сходится со скоростью геометрической прогрессии, что считается хорошим уровнем. Отметим, что в последнее время повысился интерес к этому классу задач ввиду приложений в нелинейных перепараметризованных системах глубокого обучения с перепараметризацией \cite{belkin}. Существуют и примеры седловых задач такого типа, возникающие в анализе данных (специальные робастные вариации метода наименьших квадратов) \cite{garg}. В работе рассматриваются подходы к седловым задачам с двyсторонним вариантом yсловия Поляка-Лоясиевича на базе градиентного метода с неточной информацией и предлагается правило остановки на основе малости нормы неточного градиента внешней подзадачи. Достижение этого правила в сочетании с подходящей точностью решения вспомогательной подзадачи гарантирует достижение приемлемого качества исходной седловой задачи. Обсуждаются результаты численных экспериментов для различных седловых задач для иллюстрации эффективности предложенного метода, в том числе по сравнению с доказанными оценками скорости сходимости.
\end{abstract}

\label{sec:Chapter1} \index{Chapter1}
\section{Постановка задачи}

Будем рассматривать задачи нахождения седловой точки $(x^*, y^*)$ 
\begin{align}
f(x^*, y^*) = \min\limits_{x} \max\limits_{y} f(x, y),
\end{align}
где $f$ удовлетворяет двустороннему вариантy yсловия Поляка-Лоясиевича (PL-условию), 
\begin{align}
\begin{cases}
 \|\nabla_xf(x, y)\|^2 \ge 2\mu_1(f(x,y) - \min\limits_{x} f(x,y))\\
 \|\nabla_yf(x, y)\|^2 \ge 2\mu_2(\max\limits_{y} f(x,y) - f(x,y))
\end{cases}
\end{align}
а также условиям Липшица градиента относительно евклидовой нормы.
\begin{align}\label{lipshits}
\begin{cases}
\|\nabla_x f(x_1, y) - \nabla_x f(x_2, y)\| \leq L_{11}\|x_1 - x_2\|\\
\|\nabla_x f(x, y_1) - \nabla_x f(x, y_2)\| \leq L_{12}\|y_1 - y_2\|\\
\|\nabla_y f(x, y_1) - \nabla_y f(x, y_2)\| \leq L_{22}\|x_1 - x_2\|\\
\|\nabla_y f(x_1, y) - \nabla_y f(x_2, y)\| \leq L_{12}\|x_1 - x_2\|
\end{cases}
\end{align}
для всяких $x$ и $y$.

Напомним, что для задач минимизации без ограничений 
$$\min_{x} f(x),\text{ }x \in \mathbb{R}^n$$
условие градиентного доминирования Поляка-Лоясиевича (далее yсловие PL), которое гарантирyет сходимость градиентного спуска со скоростью геометрической прогрессии для достаточно гладких задач без дополнительных предположений о выпуклости фyнкции. При этом в оценки скорости сходимости не входит параметр размерности пространства. Точнее говоря, функция $f(\cdot)$ удовлетворяет условию PL, если она имеет непустое множество решений и конечное оптимальное значение $f^*$, а также существует некоторая $\mu > 0$ такая, что 
$$\|\nabla f(x)\|^2 \geq 2\mu(f(x) - f^*)\text{, для всех } x.$$ 
Известно, что:
\begin{itemize}
    \item PL-условие является более слабым по сравнению с сильной выпуклостью, т.е. оно верно для более широкого класса функций.
    \item В частности, PL-условие верно для невыпyклых задач, связанных с перепараметризованными системами в глyбоком обyчении. \cite{belkin}
\end{itemize}

\section{Примеры задач рассматриваемого класса}
Существует некоторое количество классов задач, удовлетворяющих описанным выше условиям, однако для отдельно взятой подходящей функции вычисление параметров $L_{ij}$ и $\mu_i$ тоже может являться достаточно сложной задачей. Давайте рассмотрим несколько примеров удовлетворяющих требуемым ограничениям. 

\subsection{Пример 1}
$f(x,y)=F(Ax, By)$, где $F(\cdot,\cdot)$ сильно-выпукло-сильно-вогнутая и A, B~--- произвольные матрицы, удовлетворяет двустороннему PL-условию.

Пример взят из статьи \cite{yang}.
\subsection{Пример 2}
Невыпукло-невогнутая $f(x,y)=x^2+3\sin^2x\sin^2y-4y^2-10\sin^2y$, удовлетворяет двустороннему PL-условию с $\mu_1 = 1/16, \mu_2 = 1/14$ и условиям Липшица градиента $L_{11}=8, L_{22}=28$.

Пример взят из статьи \cite{yang}.
\subsection{Пример 3}
Задача Robust least squares (RLS) состоит в минимизации суммы квадратов расхождения между предсказанными значениями и фактическими значениями целевой переменной, но с учетом возможных выбросов в данных. Другими словами, RLS ищет оптимальную модель, которая максимально точно описывает данные, но при этом устойчива к наличию шумов и выбросов в использyемых данных. Задача RLS часто используется в машинном обучении и статистике для построения устойчивых моделей регрессии и классификации.\\

Можно рассмотреть следующий пример задачи RLS \cite{garg}: 
$$F(x,y) := \|Ax-y\|_M^2 - \lambda\|y-y_0\|_M^2\text{, где }\|\cdot\|_M^2 = x^TMx,$$ 
где M положительно полуопределена и $\lambda > 1$. В статье \cite{garg} показано, что такая задача удовлетворяет двустороннему PL-условию.

\section{Цель статьи}
Для решения таких задач часто используется градиентный метод с возможностью использования неточной информации о градиенте, который будет подробно рассмотрен в следующей главе. Теоретические оценки количества итераций, необходимых для достижения приемлемого качества решения с помощью такого градиентного метода в общем случае неулучшаемы. Однако, по-видимомy, на практике они предполагают значительное завышение затрат в сравнении с реальной необходимостью. Мы отправляемся от статьи \cite{stonyakin}, в которой уже было предложено правило ранней остановки для задач минимизации. Однако для седловых задач такого еще нет. Поэтому, основной целью данной работы стала разработка правил ранней остановки градиентных методов для седловых задач, которые могли бы гарантировать достижение приемлемого качества точки выхода по функции, а также провести эксперименты на различных седловых задачах и сравнить численные результаты предложенного метода с доказанными оценками, чтобы убедиться в эффективности предложенных методов.

\section{Градиентный метод с неточной информацией и его использование в седловых задачах}

Считаем, что значения параметров $L_{11},\;L_{12},\;L_{21},\;L_{22},\;\mu_1,\;\mu_2,\;\gamma$, которые будут задавать точность решения внутренней подзадачи, известны и положительны.\\
Будем сводить поставленную задачу к задаче минимизации вспомогательной функции вида
\begin{align}
    g(x) = \max\limits_y f(x,y).
\end{align}
Тогда согласно лемме A.5 \cite{nouiehed} получаем
\begin{align}\label{lipshits_for_outer}
    \left\|\nabla g(x_1) - \nabla g(x_2)\right\| &\leq L\left\|x_1 - x_2\right\|
\end{align}
при
\begin{align}\label{real_L}
L = L_{11} + \frac{L^2_{12}}{\mu_2}.
\end{align}

При этом обычный градиентный спуск в этом случае не применим. Это связано с тем, что на практике возможно, как правило, лишь приближённо для всякого $x$ найти значение градиента $g(x)$. Поэтому считаем, что в любой точке $x$ нам доступно значение неточного градиента $\widetilde{\nabla} g(x) = \nabla_x f(x, y)$, причем $\|\widetilde{\nabla} g(x) - g(x)\| \leq \Delta$, при некотором фиксированном положительном $\Delta=L_{12}\gamma$.\\
Тогда \eqref{lipshits_for_outer} означает, что
\begin{equation}\label{estimation_inexact_grad}
    g(x) - g^* \leq \frac{1}{\mu_1}(\|\widetilde{\nabla} g(x)\|^2 + \Delta^2) \quad \forall x \in \mathbb{R}^n,
\end{equation}
поэтому для всякого $x$ верно
$$
\|\widetilde{\nabla} g(x)\|^2 \geq \mu_1(g(x) - g^*) - \Delta^2.
$$
К задаче минимизации $g$ будем применять градиентный метод вида
\begin{align}\label{outer}
x_{k+1} &= x_{k} - \frac{1}{L} \widetilde{\nabla} g(x_k)
\end{align}
для <<внешней задачи>> и
\begin{align}\label{inner}
y_{l+1} &= y_{l} + \frac{1}{L_{22}} \nabla_y f(x_k, y_l)
\end{align}
для каждой из <<внутренних>> подзадач, необходимых для нахождения приближенного значения градиента для внешней задачи с заданной точностью решения по аргументу.

Если $\|y_k - y^*\| \leq \gamma$ и $\gamma$ достаточно мал, а также
$$
g(x_k) = \max_y f(x_k, y) = f(x_k, y_k^*),
$$
тогда имеем:
\begin{align}\label{grad}
\|\widetilde{\nabla} g(x_k) - \nabla g(x_k)\| = \|\nabla f(x_k, y_k) - \nabla_x f(x_k, y_k^*)\| \
\leq L_{12}\|y_k - y_k^{*}\| \leq L_{12}\gamma.
\end{align}

Таким образом, в \eqref{grad} $\widetilde{\nabla} g(x_k)$ представляет собой аддитивно неточный градиент $g$ в точке $x_k$ с $\Delta = L_{12}\gamma$, а $\gamma$ зависит от ошибки решения вспомогательной проблемы для $y$.\\
Лемма A3 в \cite{yang} (стр. 16) утверждает, что
$$
\left\|\nabla g(x) \right\|^2 \ge 2\mu_1 (g(x) - g(x^*)),
$$
т.е $g$ удовлетворяет обычному условию Поляка-Лоясиевича.

Получим оценку на количество итераций для <<внутренней>> подзадачи, гарантирующее приемлемую её точность.

Из известного результата сходимости со скоростью геометрической прогрессии по аргументу для L-гладких функций, удовлетворяющих PL-условию (а значит, и квадратичному росту) получаем, что
$$||y_k - y_k^*||^2 \leq \frac{2}{\mu_2}(f(x_k, y_k^*)-f(x_k, y_k)) \leq $$
$$ \leq \frac{2}{\mu_2}(1 - \frac{\mu_2}{L_{22}})^p\left(f(x_k, y_k^*) - f(x_k, y_0)\right) \leq \frac{2}{\mu_2}(1 - \frac{\mu_2}{L_{22}})^pC_2$$
для некоторого $C_2 > 0$.

Если потребовать, чтобы
$$L_{12}^2\gamma^2 \geq \frac{2}{\mu_2}(1 - \frac{\mu_2}{L_{22}})^p C_2,$$
тогда
$$(1 - \frac{\mu_2}{L_{22}})^p \leq \exp^{-\frac{\mu_2 p}{L_{22}}} \leq \frac{L_{12}^2 \gamma^2 \mu_2}{2C_2},$$
и при
\begin{align}\label{inner_upper_constraint}
p \geq \left\lceil \frac{L_{22}}{\mu_2}\log\frac{2C_2}{L_{12}^2 \gamma^2 \mu_2} \right\rceil
\end{align}
получаем требуемую погрешность по аргументу.

Продолжим для <<внешней>> подзадачи. Ввиду \eqref{lipshits_for_outer} для метода \eqref{outer} получим следующее соотношение:
$$
g(x_{k+1}) \leq g(x_k) + \langle \nabla g(x_k), x_{k+1} - x_k \rangle + \frac{L}{2}\|x_{k+1} - x_k\|^2=
$$
$$
=g(x_k) - \frac{1}{L} \langle \nabla g(x_k), \widetilde{\nabla}g(x_k) \rangle + \frac{1}{2L}\|\widetilde{\nabla}g(x_k)\|^2=
$$
$$
=g(x_k) + \frac{1}{2L} (\nabla g(x_k)^2 - 2 \langle \nabla g(x_k), \widetilde{\nabla}g(x_k) \rangle + \widetilde{\nabla}g(x_k)^2) - \frac{\|\nabla g(x_k)\|^2}{2L}=
$$
$$
=g(x_k) + \frac{1}{2L} \|\nabla g(x_k) - \widetilde{\nabla}g(k_k)\|^2 - \frac{\|\nabla g(x_k)\|^2}{2L} \leq
$$
$$
\leq g(x_k) + \frac{\Delta^2}{2L} - \frac{1}{2L}\|\nabla g(x_k)\|^2,
$$
т.~е.
\begin{align}\label{outer_decrese}
g(x_{k + 1}) - g(x_k) \leq \frac{\Delta^2}{2L} - \frac{1}{2L}\|\nabla g(x_k)\|^2.
\end{align}
Из пункта 3.4 \cite{stonyakin}:

$$g(x_{k}) - g(x^*) \leq (1 - \frac{\mu_1}{L})^{k}(g(x_0) - g(x^*)) + \frac{\Delta^2}{2\mu_1} =$$
$$= (1 - \frac{\mu_1}{L})^{k}(g(x_0) - g(x^*)) + \frac{L_{12}^2 \gamma^2}{2\mu_1} = $$
$$= (1 - \frac{\mu_1\mu_2}{L_{11}\mu_2 + L_{12}^2})^{k} (g(x_0) - g(x^*)) + \frac{L_{12}^2 \gamma^2}{2\mu_1}=$$
\begin{align}\label{proof_t3}
= (1 - \frac{\mu_1\mu_2}{L_{11}\mu_2 + L_{12}^2})^{k} C_1 + \frac{L_{12}^2 \gamma^2}{2\mu_1}
\end{align}
для некоторого $C_1 > 0$, тогда
\begin{align}
\begin{cases}
C_1(1 - \frac{\mu_1}{L})^k \leq \exp^{-\frac{\mu_1}{L}k} C_1 \leq \frac{\varepsilon}{2},\\
\frac{L_{12}^2 \gamma^2}{2\mu_1} \leq \frac{\varepsilon}{2}.
\end{cases}
\end{align}
Отсюда получаем оценку на количество итераций для внешней подзадачи:

$$\exp^{-\frac{\mu_1}{L}k} \leq \frac{\varepsilon}{2C_1},$$
$$-\frac{\mu_1}{L}k \leq -\log \frac{2C_1}{\varepsilon},$$
$$k \geq \left \lceil \frac{L}{\mu_1} \log\frac{2C_1}{\varepsilon} \right \rceil.$$
А также $L_{12}^2 \gamma^2 \leq \varepsilon \mu_1$, и тогда можно переписать оценку количества итераций для решения внутренней подзадачи:
$$p \geq \left \lceil \frac{L_{22}}{\mu_2}\log\frac{2C_2}{\varepsilon \mu_1 \mu_2} \right \rceil.$$
Итого при
$$N = k \cdot p > \left \lceil \frac{L}{\mu_1} \log\frac{2C_1}{\varepsilon} \right\rceil \cdot \left\lceil \frac{L_{22}}{\mu_2}\log\frac{2C_2}{\varepsilon \mu_1 \mu_2} \right \rceil$$
мы достигнем условия $g(x_k)-g(x^*)\leq\varepsilon$.\\

Такой подход к решению для седловых задач на базе градиентного метода с неточной информацией является известным, однако нам неизвестны случаи, в которых аккуратно были бы прописаны оценки сложности этой схемы для седловых задач с двусторонним аналогом PL-условия.

\section{Правила остановки}\label{sec:Chapter3}\index{Chapter3}

В предыдущей главе мы рассмотрели подход к решению седловых задач такого типа на базе градиентного метода с неточной информацией. В этой же будет предложен критерий ранней остановки на базе малости нормы неточного градиента внешней подзадачи, который гарантирует достижение приемлемого качества решения исходной седловой задачи. А также критерий остановки для внутренних подзадач, гарантирующий достижение требуемого качества решения по аргументу. Будет дана оценка количества итераций, необходимых для их реализации.

\subsection{Критерий остановки для внешней подзадачи}

\noindent Из \eqref{outer_decrese} ввиду
$$
\|\nabla g(x_k)\|^2 \geq \frac{\|\widetilde{\nabla}g(x_k)\|^2}{2} - L_{12}^2\gamma^2
$$
верно
$$
f(x_{k + 1}) - f(x_k) \leq \frac{L_{12}^2\gamma^2}{2L} - \frac{1}{2L}\left(\frac{\|\widetilde{\nabla}g(x_k)\|^2}{2} - L_{12}^2\gamma^2\right),
$$
откуда имеем
\begin{align}\label{outer_decrese_inexact}
g(x_{k + 1}) - g(x_k) \leq \frac{L_{12}^2\gamma^2}{L} - \frac{1}{4L}\|\widetilde{\nabla}g(x_k)\|^2.
\end{align}

Из неравенства \eqref{outer_decrese_inexact} видно, что если значение $\|\widetilde{\nabla}g(x_k)\|$ достаточно велико, то можно гарантировать, что $g(x_{k+1}) < g(x_k)$, что указывает на конечность процесса. Тем самым, для всякого $C > 2$ возникает альтернатива: или верно неравенство $\|\nabla g(x_k)\| \leq CL_{12}\gamma$, и это гарантирует достижение приемлемого качества точки выхода $x_k$ по функции в силу PL-условия, или же
$$
g(x_{k+1}) - g(x_k) < -\frac{L_{12}^2\gamma^2}{L}\left(\frac{C^2}{4} - 1\right).
$$

Тем самым, за конечное число шагов градиентного метода \eqref{outer} возможно получить $x_k$ такое, что значение $g(x_k)$ достаточно близко к минимальному $g(x^*)$. Выберем для определённости $C = \sqrt{6}$ (чтобы получить <<удобный>> коэффициент) и будем рассматривать 2 сценария:
\begin{enumerate}
\item
$\|\widetilde{\nabla}g(x_k)\| > L_{12}\gamma\sqrt{6}$, откуда с учетом \eqref{outer_decrese_inexact} получаем
\begin{align}\label{first_situation}
g(x_{k+1}) - g(x_k) < -\frac{L_{12}^2\gamma^2}{2L}.
\end{align}
\item
\begin{align}\label{criteria_2}
\|\widetilde{\nabla}g(x_k)\| \leq L_{12}\gamma\sqrt{6},
\end{align}
откуда ввиду \eqref{estimation_inexact_grad} имеем
\begin{align}\label{estimation_delta_f}
g(x_{k+1}) - g(x^*) \leq \frac{7L_{12}^2\gamma^2}{\mu_1}.
\end{align}
\end{enumerate}

Будем считать оценку \eqref{estimation_delta_f} приемлемой и договоримся обрывать процесс \eqref{outer} в случае, если выполнено \eqref{criteria_2}.

\subsection{Критерий остановки для внутренних подзадач}

Обозначим $t(y) = f(*, y)$~--- внутренняя подзадача для фиксированного $x$,

$$
L_{22}(y_{k+1}-y_k)=\Delta t(y_k),
$$
$$
t(y_{k+1}) \geq t(y_k) + \langle \Delta t(y_k), y_{k+1}-y_k \rangle - \frac{L_{22}}{2} \|y_{k+1}-y_k\|^2,
$$
$$
t(y_{k+1}) - t(y_k) \geq \frac{1}{L_{22}}\|\nabla t(y_k)\|^2 - \frac{1}{2L_{22}}\|\nabla t(y_k)\|^2,
$$
$$
t(y_{k+1}) - t(y_k) \geq \frac{1}{2L_{22}}\|\nabla t(y_k)\|^2.
$$
Пусть $\|\nabla t(y_k)\| \leq C_3$, тогда из PL-условия и условия квадратичного роста:
\begin{equation}\label{inner_c}
\frac{\mu_2}{2} \|y_k-y^*\|^2 \leq t^*-t(y_k) \leq \frac{1}{2\mu_2}C_3^2.
\end{equation}
Ввиду требуемого $\|y_k-y^*\| \leq \gamma$ возьмем $C_3 = \mu_2 \cdot \gamma$.

Таким образом, возникает альтернатива: или градиент становится достаточно мал, и ввиду PL-условия нам удается достигнуть необходимого качества решения по аргументу, либо продолжаем делать шаги градиентного метода.

Также можно использовать $\hat{y}$ с предыдущего шага как начальную точку для следующего. Тогда на практике для задач, в которых требуется большая точность, критерий остановки будет срабатывать <<почти сразу>>.

\subsection{Итоговая схема}

К седловой задаче для $f$ применяются методы вида
\begin{equation}\label{method_1}
    x_{k+1} = x_{k} - \frac{\mu_2}{L_{11}\mu_2 + L_{12}^2} \widetilde{\nabla} g(x_k) = x_{k} - \frac{\mu_2}{L_{11}\mu_2 + L_{12}^2} \nabla_x f(x_k, y_m) 
\end{equation}
для <<внешней>> задачи, и для вычисления $\widetilde{\nabla} g(x_k)$~---
\begin{equation}\label{method_2}
y_{m+1} = y_{m} + \frac{1}{L_{22}} \nabla_y f(x_k, y_m)
\end{equation}
для каждой из <<внутренних>> подзадач.\\

Для каждого из методов \eqref{method_1}, \eqref{method_2} предлагаются правила ранней остановки
\eqref{rule_1}, \eqref{rule_2} соответственно:
\begin{equation}\label{rule_1}
    \|\nabla_x f(x_k, y_m)\| \leq L_{12}\gamma \sqrt{6},
\end{equation}
\begin{equation}\label{rule_2}
    \|\nabla_y f(x_k, y_m)\| \leq \mu_2 \gamma.
\end{equation}

В таком случае справедливы следующие теоремы.

\subsubsection{Теорема 1}

\begin{theo}
Пусть на $p$-й итерации градиентного метода \eqref{method_2} впервые выполнен критерий остановки \eqref{rule_2}. Тогда для точки выхода $\hat{y} = y_p$ гарантированно будет верно неравенство
\begin{equation}
    \|y^*(x_k)-\hat{y}\| \leq \gamma.
\end{equation}
При этом справедлива следующая оценка количества итераций до полной остановки:
\begin{equation}
p \leq \left \lceil \frac{L_{22}}{\mu_2}\log\left(\frac{2C_2}{L_{12}^2\gamma^2\mu_2^2}\right) \right \rceil.
\end{equation}
\noindent\rule{\textwidth}{1pt}
$\forall x \mapsto f(x, y^*) - f(x, y_0) \leq C_2$ для некоторого $C_2 \geq 0$.
\end{theo}

Что следует из \eqref{inner_upper_constraint}, а также предложенного критерия остановки для внутренней подзадачи.
\newpage

\subsubsection{Теорема 2}

\begin{theo}
Для градиентного метода \eqref{method_1} критерий остановки \eqref{rule_1} выполнен при
\begin{equation}
N \leq \left \lceil 2C_1\frac{L_{11} + \frac{L_{12}^2}{\mu_2}}{L_{12}^2\gamma^2} \right \rceil.
\end{equation}
\noindent\rule{\textwidth}{1pt}
$g(x_0) - g(x^*) \leq C_1$ для некоторого $C_1 \geq 0$.
\end{theo}

Пусть критерий остановки \eqref{rule_1} не выполнен для всех $k$ от $0$ до $N - 1$. Тогда $\widetilde{\nabla} g(x_k) > \sqrt{6}L_{12}\gamma$. Но тогда, просуммировав выражения \eqref{first_situation} для всех $k$, получим следующее соотношение:
$$
g(x_0)-g(x^*) \geq g(x_0)-g(x_N) = \sum\limits_{k=0}^{N-1}(g(x_k)-g(x_{k+1}))>\frac{NL_{12}^2\gamma^2}{2L}.
$$

Подставив $L$ из \eqref{real_L}, получаем требуемое соотношение.

Эта оценка является сильно завышенной, однако ее удобно использовать при отсутствии информации о параметре $\mu_1$.

\subsubsection{Теорема 3}

\begin{theo}
Пусть градиентный метод \eqref{method_1} работает либо
\begin{equation}
N_* = \left \lceil \frac{L_{11}\mu_2 + L_{12}^2}{\mu_1\mu_2} \log\left(\frac{C_1\mu_1} {6L_{12}^2\gamma^2}\right) \right \rceil
\end{equation}
шагов, либо при некотором $N \leq N_*$ на  $N$-й итерации метода \eqref{method_1} впервые выполнен критерий остановки \eqref{rule_1}. Тогда для точки выхода $\hat{x} = x_N$ гарантированно будет верно неравенство
\begin{equation}
\|\hat{x} - x^*\| \leq \frac{L_{12}\gamma \sqrt{14}}{\mu_1}.
\end{equation}
\noindent\rule{\textwidth}{1pt}
Для некоторых $C_1, C_2$.
\end{theo}

Пусть также критерий остановки \eqref{rule_1} не выполнен для всех $k$ от $0$ до $N - 1$ и $\widetilde{\nabla} g(x_k) > \sqrt{6}L_{12}\gamma$. Тогда ввиду \eqref{estimation_delta_f} и \eqref{proof_t3} достаточно потребовать, чтобы
$$
\left(1 - \frac{\mu_1\mu_2}{L_{11}\mu_2 + L_{12}^2}\right)^{k} C_1 \leq 6\frac{L_{12}^2 \gamma^2}{\mu_1},
$$
откуда получаем требуемую оценку на $N_{*}$. Далее, при выполненном критерии верно \eqref{estimation_delta_f}. И из условия квадратичного роста:
$$
\|x_N-x_{*}\|^2\leq \frac{2}{\mu_1}(g(x_N)-g(x^*))\leq 14\frac{L_{12}^2\gamma^2 }{\mu_1^2}.
$$
Извлекая квадратный корень из выражений, получим требуемое.

\section{Вычислительные эксперименты}

\noindent Будем проводить эксперимент на следующем примере:
\begin{equation}\label{example_1}
    f(x,y)=x_1^2+x_2^2+x_3^2+3\sin^2x_1\sin^2y_1-4y_1^2-3y_2^2-2y_3^2-10\sin^2y_1
\end{equation}
Функция \eqref{example_1} удовлетворяет необходимым условиям с коэффициентами
\begin{itemize}
    \item $L_{11} = 8$,
    \item $L_{22} = 28$,
    \item $\mu_1 = 1/16$,
    \item $\mu_2 = 1/14$.
\end{itemize}
Константы $L_{12},\;L_{21} \leq L_{22} = 28$. Также зафиксируем константы $C_1=C_2=100$. Для данной задачи оптимальным решением является точка $\Vec{0}$.

\begin{center}
Результаты при случайном выборе $y$.
\begin{tabular}{ ||c|c|c|c|c||c|c|| }
 \hline
 $\gamma$ & $f(\hat{x}, \hat{y})-f^*$ & $N$ & $\sum p_k$ & $avg_p$ & $N^*$ & $p$ \\
 \hline
 $10^{-3}$ & $1.1749 \cdot 10^{-3}$ & 40921 & 3921481 & 95.83 & 1751204 & 6950 \\
 \hline
 $10^{-5}$ & $1.1752 \cdot 10^{-7}$ & 69525 & 8737867 & 125.67 & 3369866 & 10560 \\
 \hline
 $10^{-8}$ & $1.1755 \cdot 10^{-13}$ & 103487 & 17644048 & 170.49 & 5797859 & 15976 \\
 \hline
\end{tabular}
\end{center}
\begin{center}
Результаты при выборе $y_0 = y_{opt}$ на предыдущем шаге.
\begin{tabular}{ ||c|c|c|c|c||c|c|| }
 \hline
 $\gamma$ & $f(\hat{x}, \hat{y})-f^*$ & $N$ & $\sum p_k$ & $avg_p$ & $N^*$ & $p$ \\
 \hline
 $10^{-3}$ & $1.1752 \cdot 10^{-3}$ & 45266 & 80367 & 1.77 & 1751204 & 6950 \\
 \hline
 $10^{-5}$ & $1.1753 \cdot 10^{-7}$ & 65430 & 129758 & 1.98 & 3369866 & 10560 \\
 \hline
 $10^{-8}$ & $1.1755 \cdot 10^{-13}$ & 106811 & 237818 & 2.23 & 5797859 & 15976 \\
 \hline
\end{tabular}
\end{center}

Здесь $p_k$~--- это количество шагов градиентного метода внутренней подзадачи на $k$-ой итерации внешней, а $avg_p$~--- среднее шагов внутренних на одной итерации внешней.

По итогам экспериментов видно, что градиентный метод с критериями остановки намного эффективнее, чем выполнение теоретически необходимого количества итераций. Если считать количество суммарных операций взятия градиента самой ресурсоемкой частью алгоритма, то теоретически потребовалось бы провести $N^* \cdot (p + 1)$ таких операций, тогда как на практике требуемая точность достигается при $N + \sum p_k$, что более чем в $10^4$ раз меньше в случае со случайным $y_0$ и в $10^5$ раз~--- при использовании информации с предыдущего запуска.

\end{document}